\documentclass[10 pt]{article}
\usepackage{geometry}
\usepackage{amsthm}
\usepackage{amsmath}
\usepackage{amssymb}
\usepackage{bm}

\newtheorem{theorem}{Theorem}[section]

\newtheorem{proposition}[theorem]{Proposition}

\begin{document}

\title{\textbf{Characterizing Logarithmic Bregman Functions \footnote{Submitted in \textit{Stat}}}}
\author{Ray, S.$^{1}$, Pal, S.$^{2}$, Kar, S.~K.$^{3}$ and Basu, A.$^{4}$\\$^{1}${\small Stanford University}, $^{2}${\small Iowa State University}\\$^{3}${\small University of North Carolina, Chapel Hill}, $^{4}${\small Indian Statistical Institute}}  
\date{}
\maketitle

\begin{abstract}
Minimum divergence procedures based on the density power divergence and the logarithmic density power divergence have been extremely popular and successful in generating inference procedures which combine a high degree of model efficiency with strong outlier stability. Such procedures are always preferable in practical situations over procedures which achieve their robustness at a major cost of efficiency or are highly efficient but have poor robustness properties. The density power divergence (DPD) family of \cite{B01} and the logarithmic density power divergence (LDPD) family of \cite{J02} provide flexible classes of divergences where the adjustment between efficiency and robustness is controlled by a single, real, non-negative parameter. The usefulness of these two families of divergences in statistical inference makes it meaningful to search for other related families of divergences in the same spirit. The DPD family is a member of the class of Bregman divergences, and the LDPD family is obtained by log transformations of the different segments of the divergences within the DPD family. Both the DPD and LDPD families lead to the Kullback-Leibler divergence in the limiting case as the tuning parameter $\alpha \rightarrow 0$. In this paper we study this relation in detail, and demonstrate that such log transformations can only be meaningful in the context of the DPD (or the convex generating function of  the DPD) within the general fold of Bregman divergences,  giving us a limit to the extent to which the search for useful divergences could be successful.
\end{abstract}









\section{Introduction}\label{sec1}

Statistical distances, or divergences, between two probability density functions are extensively used in many branches of 
science including statistics, machine learning and mathematical information theory. 
In this paper we will deal with the interdependence between some useful divergences and their 
characterizing properties (rather than their utility in any particular field of application). We believe that scientists in many different 
areas of application (including those mentioned above) will find these characterizations to be useful.

We will, in particular, be interested in two interlinked classes of divergences which have been extensively used in many different scientific domains. 
These two families of divergences are the DPD and LDPD classes (or the BHHJ and the JHHB classes, or the type 1 and type 0 classes, or 
the $\beta$-divergences and the $\gamma$-divergences) respectively; see, for example, \cite{J02},~\cite{F01},~\cite{C01}
We will define these two classes of divergences in Section~\ref{sec2}. Although the forms 
are very similar and the divergences have several obvious links, these divergences have been arrived at through somewhat 
different considerations; see \cite{B02}. However the apparent similarity in form prompts one to explore whether there exist other 
pairs of divergences  with similar links. 

The DPD family is a member of the class of Bregman divergences; the LDPD is not. The DPD is a single-integral, non-kernel divergence in 
the nomenclature of \cite{J01}. The LDPD is not a single-integral but is a non-kernel divergence. The non-kernel divergences 
have also been called \textit{decomposable} in the literature~\cite{B03}. The divergences within the DPD family have been shown to possess 
strong robustness properties in statistical applications. The LDPD family is also a very useful one in statistical inference. The 
applicability of the latter class of divergences in mathematical information theory has been explored by \cite{K02},~\cite{K03}. \cite{F01} have also stressed upon the importance of the generalized cross entropy functions corresponding to the LDPD, as well as 
possible advantages of the LDPD over the DPD; see also \cite{K01} for some counter-views in this connection. In general 
both families have important statistical uses.

A Bregman divergence is generated by a suitable convex function $B(\cdot)$. For the density power divergence 
\begin{align}
\label{defineDPD}
B(y) := \frac{y^{1+\alpha}-y}{\alpha}, \; y>0,
\end{align}
is the relevant convex function and $\alpha$ is a non-negative tuning parameter.  The LDPD is obtained by piecewise application of  the logarithmic 
transformation on different segments of the divergences from the DPD family.  However the origin of the LDPD in the spirit of the weighted method
of moments can also be considered as a follow-up of the work of \cite{W01}; see \cite{J02},~\cite{B02} for a description of this development. A recent work in~\cite{R01} has characterized the class of all possible functions which can be applied piecewise on different segments of the divergences from the DPD family to generate another valid divergence, logarithmic function being one particular member of this class. Although the natural question of whether any other member of this class generates a divergence with nice inferential properties remains unexplored.

In this note we will explore which other members of the class of Bregman divergences (or, equivalently what other convex functions),  
produce legitimate and useful divergence measures by log transforms in the spirit of the DPD-LDPD connection. In light of the fact 
that the DPD and the LDPD both have various meaningful scientific uses, such a characterization could be of great value for both theoretical 
and applied research.

\section{Bregman Divergence}\label{sec2}

\subsection{Usual Bregman Divergence}
Let $B: [0,\infty) \longrightarrow [0,\infty) $ be a twice differentiable, strictly convex function, and 
$(\mathbb{R},\mathcal{B}_{\mathbb{R}}, \mu)$ be a measure space on the real line. We define the set
$$ \mathcal{L}_{\mu} := \left\{ f : \mathbb{R} \longrightarrow \mathbb{R}\; \big \rvert \; f \geq 0, \;f \;\text{measurable}, \;
 \int_{\mathbb{R}} f \;d\mu = 1 \right\}.$$
We want to define a divergence on $\mathcal{L}_{\mu}$, i.e. a function $D : \mathcal{L}_{\mu} \times \mathcal{L}_{\mu} 
\longrightarrow \mathbb{R} \cup \left\{ -\infty, \infty \right\} $ such that
(i) $D(g, f) \geq 0$ for all $g, f \in \mathcal{L}_{\mu},$ and (ii) $D(g, f) =0$ if and only if $g=f \; a.e.[\mu].$
Given a convex function $B(\cdot)$, one can define such a divergence, say $D_{B}$, exploiting the convexity property as follows. Note that 
the convexity of $B(\cdot)$ ensures that 
$$ B(g(x)) - B(f(x)) - \left(g(x)-f(x)\right)B^{\prime} \left(f(x)\right) \geq 0,$$
for all $x \in \mathbb{R},$
and the inequality is strict provided $g(x)\neq f(x)$. 
(Since the support of the function $B(\cdot)$ is 
$[0, \infty)$, $B^{\prime}(0)$ is to be interpreted as the right hand derivative at zero).
Therefore,
$$ \int_{\mathbb{R}} \big[ B(g)-B(f)-B^{\prime}(f)(g-f) \big] \geq 0~~{\rm for~all}~  g, f \in \mathcal{L}_{\mu}, $$
and the equality holds if and only if $g = f\; a.e.[\mu]$. For the sake of brevity we have suppressed the dummy variable $x$ in the expressions of the integrations through out this paper. All the integrations in this paper are with respect to the measure $\mu$. 
Therefore,
\begin{align}\label{bregman_divergence}
D_B(g, f) :=  \int_{\mathbb{R}} \big[ B(g)-B(f)-B^{\prime}(f)(g-f) \big]~~ \forall~ g, f \in \mathcal{L}_{\mu}
\end{align}
defines a divergence. The $B(\cdot)$ function defined in~(\ref{defineDPD}) is one of the most common examples of a suitable convex function in this context.  
 In this case the Bregman divergence simplifies to 
\begin{align}
D_B(g, f) = \frac{1}{\alpha}\int_{\mathbb{R}} \big[ g^{\alpha+1} - f^{\alpha+1} - (\alpha+1)f^{\alpha} (g-f) \big]
&= \frac{1}{\alpha} \int_{\mathbb{R}} \big[ g^{\alpha+1} + \alpha f^{\alpha+1} - (\alpha+1)f^{\alpha}g \big] \nonumber\\
&=  \int_{\mathbb{R}} \big[  f^{\alpha+1} - \left(1+\frac{1}{\alpha}\right)f^{\alpha}g + \frac{1}{\alpha}g^{\alpha+1} \big],
\label{equation1}
\end{align}  
which is the density power divergence (DPD) class of~\cite{B01}. As mentioned earlier, in the limiting case $\alpha \rightarrow 0$, we recover the Kullback-Leibler divergence $\int g \log(g/f)$.  

We will refer to any function $B(\cdot)$ which defines a Bregman divergence in the sense of Equation (\ref{bregman_divergence}) as a 
Bregman function. Essentially this refers to all twice differentiable, strictly convex functions defined on $[0, \infty)$. 

\subsection{Logarithmic Bregman Functions}
We have already noticed that for  $B(\cdot)$ as defined in Equation~(\ref{defineDPD}) the Bregman divergence turns out to be as defined in
Equation~(\ref{equation1}). Consider the modified expression,
\begin{align}{\label{eqn2}}
\tilde{D}_B(g, f)  =   \log \int_{\mathbb{R}} f^{\alpha+1}  - \left(1+\frac{1}{\alpha}\right) \log \int_{\mathbb{R}} f^{\alpha}g  +  \frac{1}{\alpha}\log \int_{\mathbb{R}} g^{\alpha+1},
\end{align}
which turns out to be non-negative whenever   
\begin{align}\label{eqn1}
\int_{\mathbb{R}} f^{\alpha}g \leq {\left( \int_{\mathbb{R}} f^{\alpha+1}\right)}^{\frac{\alpha}{1+\alpha}} {\left( \int_{\mathbb{R}} 
g^{\alpha+1}\right)}^{\frac{1}{1+\alpha}}. 
\end{align}	
The inequality in Equation~(\ref{eqn1}) is seen to hold from an application of \textit{Holder's inequality} on the functions $f^{\alpha}$ 
and $g$ with dual indices $\left(\frac{1+\alpha}{\alpha},1+\alpha\right)$. This inequality becomes an equality if and only if  
${(f^{\alpha})}^{\frac{1+\alpha}{\alpha}} = g^{\alpha+1}, \;\; a.e.[\mu]$, i.e. $f=g, \;\; a.e.[\mu]$. This confirms that $\tilde{D}_B$ 
is indeed a divergence. Once again, in the limiting case $\alpha \rightarrow 0$, we recover the Kullback-Leibler divergence $\int g \log(g/f)$. The resulting divergence family in~(\ref{eqn2}) is known as the logarithmic density power divergence (LDPD) class.

A closer look at Equation~(\ref{eqn2}) shows that in relation to the function $B(y) = (y^{1+\alpha} - y)/\alpha$, the divergence has the structure
\begin{align}\label{long_equation}                                                              
\tilde{D}_B(g,f)  &=  \alpha_0 \log \int_{\mathbb{R}} \frac{1}{\alpha_0} [B^{\prime}(f)f-B(f)+B(0)] + \alpha_2 \log \int_{\mathbb{R}} 
\frac{1}{\alpha_2} [B(g) - B(0) - B^{\prime}(0)g]  \nonumber \\
& \hspace{ 3 in } - \alpha_1 \log \int_{\mathbb{R}} \frac{1}{\alpha_1} [B^{\prime}(f)-B^{\prime}(0)]g, 
\end{align}
where $\alpha_0,\alpha_1,\alpha_2$ are some positive constants. The adjustments using $B(0),B^{\prime}(0)$ are to make sure that all the integrands 
are non-negative, as we shall see later.  For the case of $B(y)=(y^{\alpha+1} -y)/\alpha$, the choices $\alpha_0=1, \alpha_1=(1+\alpha)/\alpha, 
\alpha_2=1/\alpha$ gives us the divergence.

Our interest lies in determining whether there exists any other such strictly convex twice differentiable non-negative function $B(\cdot)$ on 
the non-negative real line for which there exists suitable choices of indices $(\alpha_0,\alpha_1,\alpha_2)$ which will make $\tilde{D}_B$, 
as defined in Equation~(\ref{long_equation}), a proper divergence. We will refer to such functions $B(\cdot)$, which admit proper divergences under the 
indicated piecewise log transformations,  as \textit{logarithmic Bregman functions},  or in short, as LBFs. Similarly, the valid divergences of the form given in Equation (\ref{long_equation}) corresponding to legitimate LBFs are denoted as {\it logarithmic Bregman divergences}.

Before we try to characterize the LBFs, let us explain the adjustments made in the terms of the right hand side in~(\ref{long_equation}).
Notice that, strict convexity of $B$ gives, 
\begin{align*}
B(y) - B(0)  =    \int_{0}^y B^{\prime}(t) \;dt \leq \int_{0}^y B^{\prime}(y) \;dt = yB^{\prime} (y), 
\end{align*}
i.e., $ yB^{\prime}(y) - B(y)+B(0) \geq 0$ for all $y \geq 0.$ Equality holds in the last equation if and only if $B^{\prime}(y)=
B^{\prime}(t)$ for all $0 \leq t \leq y$; by the strict convexity of $B(\cdot)$, this is possible if and only if $y = 0$. This 
guarantees that the first term in Equation~(\ref{long_equation}) is well-defined. Similar arguments show that $B(y)-B(0)-B^{\prime}(0)y \geq 0,$
for all $y \geq 0$, with equality if and only if $y = 0$. This fact guarantees that the second term in Equation~(\ref{long_equation}) is 
well-defined. Strict convexity also guarantees the well-definedness of the third term in Equation~(\ref{long_equation}).

\subsection{Conditions for $\tilde{D}_B$ Being a Divergence}
In the following, we shall work with $\mu$ being the Lebesgue measure on the real line. 
Let $B(\cdot)$ be a given Bregman function which generates, for given densities $g$ and $f$, the Bregman divergence $D_B(g, f)$. It is 
a straightforward matter to check that for all real constants $a,b$; the function $B_{a,b}(y):=B(y)+ay+b$ is an equivalent Bregman function in the sense that it is also strictly convex, twice differentiable and satisfies
$D_B(g,f)=D_{B_{a,b}}(g,f)$, for all densities $f,g$. Thus $B(\cdot)$ and $B_{a,b}(\cdot)$ generate the same Bregman divergence. Moreover,  we also have $\tilde{D}_B(g,f)=\tilde{D}_{B_{a,b}}(g,f)$, indicating that if $B$ is a LBF then $B_{a,b}$ is also an equivalent LBF. In particular, 
if we consider
\begin{align}\label{defstar}
 B_{std}(y) := B(y)-B(0)-B^{\prime}(0)y, \forall~ y \geq 0,
\end{align}
 then $B(\cdot)$ and $B_{std}(\cdot)$ are equivalent Bregman 
functions (and equivalent LBFs if $B$ is indeed an LBF);  in this case 
We also have $B_{std}(0)=0=B^{\prime}_{std}(0),$ and $B^{\prime}_{std}(y) = B^{\prime}(y)- B^{\prime}(0)$. 

We call a Bregman function (resp. an LBF) to be \textit{standard} if it satisfies the property $B(0)=0=B^{\prime}(0)$. From our discussion, we arrive at the conclusion that any Bregman function (resp. LBF) has a standard form. In particular, any Bregman function (resp. LBF) $B$
can be written as $B(x)=B_{std}(x)+ax+b$, for all $x \geq 0$ and for some $a,b$; where $B_{std}$ is its standard form. Therefore, in order to characterize all the LBFs, we can restrict our attention to the set of standard Bregman functions only. The following proposition states this characterization and is the the main result of this article.

\begin{proposition}
\label{prop}
Any standard LBF $B$ has the form $B(x)= Kx^{1+\alpha}$, for some $K, \alpha>0$. Moreover, the corresponding choices for the constants $(\alpha_0, \alpha_1,\alpha_2)$ in~(\ref{long_equation}) must satisfy the following relation; $\alpha_0=\alpha \alpha_2, \hspace{0.1 in} \alpha_1=(1+\alpha) \alpha_2.$
For any choice of $\alpha, \alpha_2, K>0$, the aforementioned triplet $(\alpha_0,\alpha_1,\alpha_2)$ and $B$ gives a valid logarithmic Bregamn divergence.
\end{proposition}


\begin{proof}
Strict convexity of $B(\cdot)$, implies that $B(\cdot)$ is strictly 
increasing and hence both $B(\cdot)$ and $B^{\prime}(\cdot)$ are positive on the positive real line, and $B(y) \leq y B^{\prime}(y)$ for 
all $y\geq 0$. Also note that as $B^{\prime}(\cdot)$ is strictly increasing, $\lim_{y \rightarrow \infty} B^{\prime}(y) = c\in (0,\infty]$ 
exists. On the other hand, 
\begin{align*}
B(y) = \int_{1}^y B^{\prime}(t) \;dt+ B(1) \geq  (y-1)B^{\prime}(1) + B(1),
\end{align*}
for all $y \geq 1$; which gives $ \lim_{y \rightarrow \infty} B(y) = \infty$, as $B^{\prime}(1) > 0$. Using this and \textit{L'Hospital's Rule}, we get, 
\begin{align}{\label{eqn4}}
\lim_{y \rightarrow \infty} \dfrac{B(y)}{y} = \lim_{y \rightarrow \infty} B^{\prime}(y) = c.
\end{align} 
Since $B$ is a standard LBF, Equation~(\ref{long_equation})  simplifies to 
\begin{align}{\label{formsmp}}
\tilde{D}_B(g, f) = \alpha_0 \log \int_{\mathbb{R}} \frac{1}{\alpha_0} [B^{\prime}(f)f-B(f)]  +\alpha_2 \log \int_{\mathbb{R}} 
\frac{1}{\alpha_2} B(g)  - \alpha_1 \log \int_{\mathbb{R}} \frac{1}{\alpha_1} B^{\prime}(f)g .
\end{align}
For this to be a proper divergence, we must have 
\begin{align*}
\alpha_0 \log \int_{\mathbb{R}} \frac{1}{\alpha_0} [B^{\prime}(f)f-B(f)] 
+ \alpha_2 \log \int_{\mathbb{R}} \dfrac{1}{\alpha_2} B(g) 
-\alpha_1 \log \int_{\mathbb{R}} \dfrac{1}{\alpha_1} B^{\prime}(f)g \geq 0, 
\end{align*}
${\rm for~ all}~ g, f,$ i.e.
\begin{align}{\label{eqn5}}
{\left(\int_{\mathbb{R}} [B^{\prime}(f)f-B(f)]\right)}^{\alpha_0}  {\left(\int_{\mathbb{R}} B(g)\right) }^{\alpha_2} \geq C \;
{\left(\int_{\mathbb{R}} B^{\prime}(f)g\right)}^{\alpha_1}, 
\end{align}
${\rm for~ all}~ g, f;$ where $C:= \dfrac{\alpha_0^{\alpha_0}\alpha_2^{\alpha_2}}{\alpha_1^{\alpha_1}}$; equality should hold if and only if $f=g \; a.e.[x]$. 

Let $U(a, b)$ represent the uniform distribution on $(a, b)$. Fix $\theta \in \mathbb{R}^{+}$,  take $f = g$ to be the 
$U\left(0, \frac{1}{\theta}\right)$ density given by $f(x) = \theta I_A(x)$, where $A = (0, 1/\theta)$ and $I_A(x)$ is the indicator function of the set $A$. 
Then,
$$ \int_{\mathbb{R}} [B^{\prime}(f)f-B(f)]= \dfrac{1}{\theta} (\theta B^{\prime}(\theta)-B(\theta)) = B^{\prime}(\theta) - \dfrac{B(\theta)}{\theta},\; \int_{\mathbb{R}} B(g) = \dfrac{1}{\theta}B(\theta) \; \text{ and } \; \int_{\mathbb{R}} B^{\prime}(f)g = \dfrac{1}{\theta}\theta B^{\prime}(\theta) =B^{\prime}(\theta).$$
Since under $g = f$ we should have equality in Equation (\ref{eqn5}), we have 
\begin{align}\label{eqn6}
\left(B^{\prime}(\theta) - \dfrac{B(\theta)}{\theta}\right)^{\alpha_0} \left(\dfrac{1}{\theta}B(\theta)\right)^{\alpha_2} 
= C (B^{\prime}(\theta))^{\alpha_1}, ~{\forall}~ \theta \in \mathbb{R}^{+},
\end{align}
or in other words,
\begin{align}{\label{eqnt}}
(\theta B^{\prime}(\theta) - B(\theta))^{\alpha_0} (B(\theta))^{\alpha_2} = C \theta^{\beta}(\theta B^{\prime}(\theta))^{\alpha_1}, ~{\forall}~  \theta \in \mathbb{R}^{+},
\end{align} 
where $\beta = \alpha_0+\alpha_2-\alpha_1$. Define $u : [0,1] \longrightarrow \mathbb{R}$ as 
 $u(x) = (1-x)^{\alpha_0}x^{\alpha_2} ~~{\rm for ~ all}~ x \in [0,1]. $
This function has a global maximum at $x=\dfrac{\alpha_2}{\alpha_0+\alpha_2}$, and the maximum value is $C_0 := \dfrac{\alpha_0^{\alpha_0} \alpha_2^{\alpha_2}}{(\alpha_0+\alpha_2)^{\alpha_0+\alpha_2}}.$ We have $ B(\theta), B^{\prime}(\theta) > 0$, for $\theta > 0$ and
$B(\theta) \leq \theta B^{\prime}(\theta), ~{\rm for~all}~   \theta >0.$ This gives $ \dfrac{B(\theta)}{\theta B^{\prime}(\theta)} \in [0,1], ~{\rm for~all}~   \theta >0.$  Equation~(\ref{eqnt}) now  says that,
\begin{align}{\label{eqnt1}}
 C_0 \geq u\left(\dfrac{B(\theta)}{\theta B^{\prime}(\theta)}\right) = C \theta^\beta \dfrac{(\theta B^{\prime} (\theta))^{\alpha_1}}{(\theta B^{\prime} (\theta))^{\alpha_0+\alpha_2}} = C(B^{\prime}(\theta))^{-\beta}, \text{ for all } \theta \in \mathbb{R}^{+}.
 \end{align}


Now consider Equation (\ref{eqn4}).  Suppose $c < \infty$. Then ~(\ref{eqn4}) gives 
$$ \lim_{y \rightarrow \infty} \left(\dfrac{B(y)}{y} - B^{\prime}(y)\right) = 0. $$  
Now take $\theta \rightarrow \infty$ in Equation (\ref{eqn6}). This gives $C c^{\alpha_1} = 0,$ i.e., $c=0$; but $c >0$, which gives a contradiction. Therefore $c= \infty $. In summary, $ \lim_{y \rightarrow 0} B^{\prime}(y) = B^{\prime}(0) = 0$ and 
$\lim_{y \rightarrow \infty} B^{\prime}(y) = \infty$.

Now get back to  (\ref{eqnt1}). If $\beta < 0$, take $\theta \rightarrow \infty$ to get $C_0 =\infty$, a contradiction. If $\beta > 0$, take take $\theta \rightarrow 0$ to get $C_0 =\infty$, a contradiction. Therefore, $\beta =0$, i.e., $\alpha_0+\alpha_2=\alpha_1$ and $C=C_0$.  (\ref{eqnt1}) now becomes,
\begin{align}{\label{eqnt2}}
u\left(\dfrac{B(\theta)}{\theta B^{\prime}(\theta)}\right) = C_0, \;\; {\rm for~all}~ \theta \in \mathbb{R}^{+}.
\end{align} 
But $C_0$ is the maximum value of $u$ attained only at the point $x=\dfrac{\alpha_2}{\alpha_0+\alpha_2}$. Therefore,
$$ \dfrac{B(\theta)}{\theta B^{\prime}(\theta)} = \dfrac{\alpha_2}{\alpha_0+\alpha_2} =: \gamma, \;\; {\rm for~all}~ \theta > 0 .$$
We get the differential equation for $B$ as $B(\theta) = \gamma\theta B^{\prime}(\theta)$.The general solution for this equation looks like
$$ B(\theta) = K \theta ^{1/\gamma}, \;\; {\rm for~all}~\theta \geq 0.$$
Note that, by definition, $0 < \gamma < 1$ and therefore $\frac{1}{\gamma} >1$. This indicates that, the only LBFs $B$ with $B(0)=B^{\prime}(0)=0$ are given by,
$$ B(x) = Kx^{1+\alpha},\;\; \text{for some} \; K,\alpha >0.$$ The constraints on the constants $(\alpha_0, \alpha_1, \alpha_2)$ follows trivially from the observations that $\alpha_0+\alpha_2=\alpha_1$ and $\alpha_2/(\alpha_0+\alpha_2) = 1/(1+\alpha).$

Plugging in $B(x)=Kx^{1+\alpha}$ and $\alpha_0=\alpha \alpha_2, \alpha_1=(1+\alpha)\alpha_2$ in~(\ref{long_equation}), and doing some algebra we obtain
\begin{align}{\label{eqn10}}
\tilde{D}_B(g, f)  =   \alpha \alpha_2\log \int_{\mathbb{R}} f^{\alpha+1}  - (1+\alpha)\alpha_2 \log \int_{\mathbb{R}} f^{\alpha}g +  \alpha_2\log \int_{\mathbb{R}} g^{\alpha+1},
\end{align}
which is nothing but a positive multiple of some member in the LDPD family, introduced in~(\ref{eqn2}); and hence a valid divergence.
\end{proof}

Recalling the transformation we made to get the standard form, we conclude that all the LBFs are given by
\begin{align}{\label{char}}
B(x) = K_1x^{1+\alpha} + K_2x + K_3,\;\; K_1, \alpha >0, K_2,K_3 \in \mathbb{R}.
\end{align}
But as we have shown in the proof of Proposition~(\ref{prop}), this $B$ and its corresponding standard form, give rise to the same divergence, which is a member of the LDPD family.

\section{Conclusions}\label{sec3}
So, in conclusion, Equation (\ref{char}) describes all the logarithmic Bregman functions, and the class of logarithmic Bregman divergences is made up only of divergences within the LDPD family and their positive multiples. Thus the DPD-LDPD pair provides the only example where the relevant logarithmic transformation starting from the Bregman divergence provides a new valid divergence.

\bigskip 


\end{document}